# INFERENCE FOR COVARIATE ADJUSTED REGRESSION VIA VARYING COEFFICIENT MODELS[1]

By Damla Şentürk and Hans-Georg Müller

*Pennsylvania State University and University of California, Davis*

We consider covariate adjusted regression (CAR), a regression method for situations where predictors and response are observed after being distorted by a multiplicative factor. The distorting factors are unknown functions of an observable covariate, where one specific distorting function is associated with each predictor or response. The dependence of both response and predictors on the same confounding covariate may alter the underlying regression relation between undistorted but unobserved predictors and response. We consider a class of highly flexible adjustment methods for parameter estimation in the underlying regression model, which is the model of interest. Asymptotic normality of the estimates is obtained by establishing a connection to varying coefficient models. These distribution results combined with proposed consistent estimates of the asymptotic variance are used for the construction of asymptotic confidence intervals for the regression coefficients. The proposed approach is illustrated with data on serum creatinine, and finite sample properties of the proposed procedures are investigated through a simulation study.

**1. Introduction.** For many statistical applications, a multiple linear regression model is a standard tool,

$$
Y_{ni} = \gamma_0 + \sum_{r=1}^{p} \gamma_r X_{nri} + e_{ni}, \tag{1}
$$

for data $(X_{nri}, Y_{ni})$, $i = 1, \ldots, n$, $r = 1, \ldots, p$, where $\gamma_0$ and $\gamma_r$ are unknown parameters, $Y_{ni}$ is the response, $X_{nri}$ is the $r$th predictor and $e_{ni}$ is the error term for the $i$th subject in the sample. An implicit assumption is that predictors and response are directly observable. However, in some situations

Received June 2003; revised April 2005.
[1]Supported by NSF Grant DMS-02-04869.

both response and predictor variables may be distorted under the influence of a confounding variable. In this paper we consider a variant of (1), where one observes contaminated versions of predictors and response. Contamination of the variables in the regression model occurs through a multiplicative factor that is determined by the value of an unknown function of an observable covariate $U$. That is, instead of observing $X_{nri}$ and $Y_{ni}$, one actually observes distorted variables $\tilde{X}_{nri}$ and $\tilde{Y}_{ni}$,

$$(2) \qquad \tilde{X}_{nri} = \phi_r(U_{ni})X_{nri}, \qquad r = 1, \ldots, p, \qquad \tilde{Y}_{ni} = \psi(U_{ni})Y_{ni}.$$

Here $\psi(\cdot)$ and $\phi_r(\cdot)$ are unknown smooth functions of the contaminating covariate $U$, and the available observations are $(U_{ni}, \tilde{X}_{nri}, \tilde{Y}_{ni})$.

An example where a model of this type is relevant are the creatinine data that are explored further in Section 5. Here serum creatinine levels are regressed on cholesterol level and serum albumin. The observed response and the two predictors are known to depend on body mass index, defined as $\text{Kg/m}^2$, which thus has a confounding effect on the regression relation. Therefore, we investigate the application of a multiplicative confounding via model (2), where the confounding variable $U$ is taken to be body mass index. "Normalization" by weight or body mass index is common in the analysis of medical data, and this refers to simply dividing the measured quantities by these confounding variables. This type of normalization implicitly assumes that the confounding is of a multiplicative nature. The adjustment considered in this paper applies to a class of more general multiplicative confounding where the effects of the confounder are modeled by unknown distorting functions $\psi(\cdot)$ and $\phi_r(\cdot)$. This leads to flexible models that include a large class of confounding mechanisms. Reasonable identifiability conditions for these functions are

$$(3) \qquad E\{\psi(U)\} = 1, \qquad E\{\phi_r(U)\} = 1, \qquad r = 1, \ldots, p,$$

corresponding to the assumption that the mean distorting effect vanishes. Additional basic assumptions are that the $(X_r, U, e)$ are mutually independent for $r = 1, \ldots, p$, and that observations made on different subjects are independent, with $E(e_{ni}) = 0$, and $\text{var}(e_{ni}) = \sigma^2$. The assumption that the underlying predictors, $X_r$, and response, $Y$, are independent of the contaminating variable $U$ is an assumption defining the proposed contamination setting through defining these unobserved, underlying variables; and for that matter it is not one that can be checked in practice. Thus, the question to be answered in practice is whether or not these independence conditions help define interpretable latent variables of interest from their observable counterparts. In our creatinine example, the latent variables are defined to be body mass index adjusted serum protein levels and cholesterol level, which are commonly used in medical studies.

The contamination of the predictor and response in a multiplicative fashion as given in (2) can alter the regression relation between the original response and predictors completely. It has also been shown for the case of simple linear regression that standard adjustment methods such as nonparametric partial regression or partial regression cannot adjust for the multiplicative contamination [11]. Therefore, a modified parameter estimation procedure is necessary, one which accounts for the multiplicative confounding effect of $U$. Such a procedure was proposed in [11], where consistent parameter estimation in the model (1)–(3) was established. This estimation procedure relies on the fact that regressing $\tilde{Y}$ on $\tilde{X}_1, \ldots, \tilde{X}_p$ gives rise to a varying coefficient model. Furthermore, a main attraction of this estimation procedure is that under the identifiability conditions of vanishing mean distorting effects, it also works for the case of additive contamination, that is, $\tilde{X}_{nri} = \phi_r(U_{ni}) + X_{nri}$, $\tilde{Y}_{ni} = \psi(U_{ni}) + Y_{ni}$, and for no contamination, that is, $\phi_r(U_{ni}) = \psi(U_{ni}) = 1$ for $r = 1, \ldots, p$. Thus, the proposed estimation procedure provides a flexible and general tool for adjustment, where the specific nature of the contamination of the variables or even its mere existence need not be known.

The aim of this paper is to derive the asymptotic distribution of these parameter estimates, and to discuss applications to confidence intervals. We show that our proposed parameter estimates are asymptotically normal, and combining this result with consistent estimation of the asymptotic variance leads to asymptotic inference.

The paper is organized as follows. In Section 2 we describe the model in detail. In Section 3 issues of estimation are discussed and the results on asymptotic inference are presented. Consistent estimates for the asymptotic variance are derived in Section 4. Applications of the proposed method to creatinine data and simulation studies are in Section 5. The proofs of the main results are assembled in Section 6, followed by the Appendix with some additional technical conditions and auxiliary results.

**2. Covariate adjustment via varying coefficient regression.** Consider the model (1)–(3). Writing $\tilde{X}_{ni} = (\tilde{X}_{n1i}, \ldots, \tilde{X}_{npi})$, the regression of the observed response on the observed predictors leads to

$$E(\tilde{Y}_{ni}|\tilde{X}_{ni}^T, U_{ni})$$
$$= E\{Y_{ni}\psi(U_{ni})|\phi_1(U_{ni})X_{n1i}, \ldots, \phi_p(U_{ni})X_{npi}, U_{ni}\}$$
$$= \psi(U_{ni})E\Big\{\gamma_0 + \sum \gamma_r X_{nri} + e_{ni}|\phi_1(U_{ni})X_{n1i}, \ldots, \phi_p(U_{ni})X_{npi}, U_{ni}\Big\}.$$

Assuming that $E(e_{ni}) = 0$ and that $(e, U, X_r)$ are mutually independent for $r = 1, \ldots, p$, the model reduces to

$$E(\tilde{Y}_{ni}|\tilde{X}_{ni}^T, U_{ni}) = \psi(U_{ni})\gamma_0 + \psi(U_{ni})\sum \gamma_r \frac{\phi_r(U_{ni})X_{nri}}{\phi_r(U_{ni})}$$

(4)
$$= \beta_0(U_{ni}) + \sum \beta_r(U_{ni})\tilde{X}_{nri}.$$

Defining the functions

(5) $$\beta_0(u) = \psi(u)\gamma_0, \qquad \beta_r(u) = \gamma_r \frac{\psi(u)}{\phi_r(u)},$$

we obtain
$$\tilde{Y}_{ni} = \beta_0(U_{ni}) + \sum \beta_r(U_{ni})\tilde{X}_{nri} + \varepsilon(U_{ni}),$$

where $\varepsilon(u) = \psi(u)e$.

We find that this is a multiple varying coefficient model, that is, an extension of regression and generalized regression models where the coefficients are allowed to vary as a smooth function of a third variable [5]. A unique feature is that both the response and predictors depend on the covariate $U$.

For varying coefficient models, Hoover, Rice, Wu and Yang [6] have proposed smoothing methods based on local least squares and smoothing splines, and recent approaches include a componentwise kernel method [13], a componentwise spline method [2] and a method based on local maximum likelihood estimates [1]. Wu and Yu [14] provide a review of recent developments. We derive asymptotic distributions for an estimation method that is tailored to this special model.

**3. Estimation and asymptotic distributions.** The estimates of the regression coefficients $\gamma_r$ will be obtained by targeting weighted averages of the smooth varying coefficient functions. Even though various smoothing methods have been proposed in the literature for the estimation of these smooth varying coefficient functions, we propose a smoothing method based on binning. The main reason for the use of the binning approach is its simplicity in targeting the desired weighted averages, rather than its performance on estimating the varying coefficient functions themselves. Nevertheless, the proposed binning approach has similarities with earlier developments for longitudinal data in Fan and Zhang [3], who use the data collected at each fixed time point to fit a linear regression, obtaining the raw estimators for the smooth varying coefficient functions.

Generalizing this idea to our independent and identically distributed data scheme, we assume that the covariate $U$ is bounded below and above, $-\infty < a \leq U \leq b < \infty$ for real numbers $a < b$, and divide the interval $[a,b]$ into $m$ equidistant intervals denoted by $B_{n1}, \ldots, B_{nm}$, referred to as bins. Given $m$, the $B_{nj}$, $j = 1, \ldots, m$, are fixed, but the number of $U_{ni}$'s falling into $B_{nj}$ is random and is denoted by $L_{nj}$. For every $U_{ni}$ falling in the $j$th bin, that is, $U_{ni} \in B_{nj}$, the corresponding observed predictors are $\tilde{X}_{n1i}, \ldots, \tilde{X}_{npi}$ and the response is $\tilde{Y}_{ni}$.

After binning the data, we fit a linear regression of $\tilde{Y}_{ni}$ on $\tilde{X}_{n1i}, \ldots, \tilde{X}_{npi}$ fusing the data falling within each bin $B_{nj}$, $j = 1, \ldots, m$. The least squares estimates of the resulting multiple regression for the data in the $j$th bin are denoted by $\hat{\beta}_{nj}^T = (\hat{\beta}_{n0j}, \ldots, \hat{\beta}_{npj})^T$. The estimators of $\gamma_{n0}$ and $\gamma_{nr}$, for $r = 1, \ldots, p$, are then obtained as weighted averages of the $\hat{\beta}_{nj}$'s, weighted according to the number of data $L_{nj}$ in the $j$th bin,

$$\hat{\gamma}_{n0} = \sum_{j=1}^{m} \frac{L_{nj}}{n} \hat{\beta}_{n0j} \tag{6}$$

and

$$\hat{\gamma}_{nr} = \frac{1}{\bar{\tilde{X}}_{nr}} \sum_{j=1}^{m} \frac{L_{nj}}{n} \hat{\beta}_{nrj} \bar{\tilde{X}}'_{nrj}, \tag{7}$$

where $\bar{\tilde{X}}_{nr} = n^{-1} \sum_{i=1}^{n} \tilde{X}_{nri}$ and $\bar{\tilde{X}}'_{nrj}$ is the average of the $\tilde{X}_{nri}$ falling in $B_{nj}$, that is, $L_{nj}^{-1} \sum_{i=1}^{n} \tilde{X}_{nri} \mathbf{1}_{\{U_{ni} \in B_{nj}\}}$ [11]. These estimates are motivated by $E\{\beta_0(U)\} = \gamma_0$ and $E\{\beta_r(U)\tilde{X}_r\} = \gamma_r E(\tilde{X}_r)$ [see (5) and (3)].

We present the asymptotic distribution of estimates $\hat{\gamma}_{n0}$ in (6), $\hat{\gamma}_{nr}$ in (7) for $\gamma_0$, $\gamma_r$ in model (1), when the number of subjects $n$ tends to infinity. As in typical smoothing applications, the number of bins $m = m(n)$ is required to satisfy $m \to \infty$, $n/(m \log n) \to \infty$ and $m/\sqrt{n} \to \infty$ as $n \to \infty$. We denote convergence in distribution by $\xrightarrow{\mathcal{D}}$ and convergence in probability by $\xrightarrow{p}$.

THEOREM 1. *Under the technical conditions* (C1)–(C7) *in Section* 6, *on event* $E_n$ [*defined in* (12)] *with* $P(E_n) \to 1$ *as* $n \to \infty$,

$$\sqrt{n}(\hat{\gamma}_{nr} - \gamma_r) \xrightarrow{\mathcal{D}} \mathbb{N}(0, \sigma_r^2), \qquad 0 \leq r \leq p,$$

*where*

$$\sigma_0^2 = \gamma_0^2 \operatorname{var}\{\psi(U)\} + \sigma^2 (\mathcal{X}^{-1})_{11} E\{\psi^2(U)\},$$

$$\sigma_r^2 = \frac{\gamma_r^2 [E(X_r^2) E\{\psi^2(U)\} - \{E(X_r)\}^2] + \sigma^2 \{E(X_r)\}^2 E\{\psi^2(U)\} (\mathcal{X}^{-1})_{rr}}{\{E(X_r)\}^2}$$

$$- \frac{2\gamma_r^2 [E\{\phi_r(U)\psi(U)\} E(X_r^2) - \{E(X_r)\}^2] + \gamma_r^2 \operatorname{var}(\tilde{X}_r)}{\{E(X_r)\}^2}, \qquad 1 \leq r \leq p,$$

*and*

$$\mathcal{X} = \begin{bmatrix} 1 & E(X_1) & \ldots & E(X_p) \\ E(X_1) & E(X_1^2) & \ldots & E(X_1 X_p) \\ \vdots & & \ddots & \vdots \\ E(X_p) & E(X_1 X_p) & \ldots & E(X_p^2) \end{bmatrix} \tag{8}$$

*is assumed to be nonsingular, according to condition* (C5) *in Section* 6.

**4. Estimating the asymptotic variance.** The observable data is of the form $(U_{ni}, \tilde{X}_{ni}^T, \tilde{Y}_{ni}), i = 1, \ldots, n$, for a sample of size $n$, where $\tilde{X}_{ni} = (\tilde{X}_{n1i}, \ldots, \tilde{X}_{npi})$ are the $p$-dimensional predictors. Correspondingly, the underlying unobservable predictors, responses and errors are $(X_{ni}^T, Y_{ni}, e_{ni}), i = 1, \ldots, n$, where $X_{ni} = (X_{n1i}, \ldots, X_{npi})$. Let $\{(U'_{njk}, \tilde{X}'_{nrjk}, \tilde{Y}'_{njk}, X'_{nrjk}, Y'_{njk}, e'_{njk}), k = 1, \ldots, L_{nj}, r = 1, \ldots, p\} = \{(U_{ni}, \tilde{X}_{nri}, \tilde{Y}_{ni}, X_{nri}, Y_{nri}, e_{nri}), i = 1, \ldots, n, r = 1, \ldots, p : U_{ni} \in B_{nj}\}$ denote the data for which $U_{ni} \in B_{nj}$, where we refer to $(U'_{njk}, \tilde{X}'_{nrjk}, \tilde{Y}'_{njk}, X'_{nrjk}, Y'_{njk}, e'_{njk})$ as the $k$th element in bin $B_{nj}$. Further let $(U'^T_{nj}, \tilde{X}'_{nj}, \tilde{Y}'^T_{nj}, X'_{nj}, Y'^T_{nj}, e'^T_{nj})$ be the data matrix belonging to the $j$th bin, where $U'_{nj} = (U'_{nj1}, \ldots, U'_{njL_{nj}})$, $\tilde{Y}'_{nj} = (\tilde{Y}'_{nj1}, \ldots, \tilde{Y}'_{njL_{nj}})$, $Y'_{nj} = (Y'_{nj1}, \ldots, Y'_{njL_{nj}})$, $e'_{nj} = (e'_{nj1}, \ldots, e'_{njL_{nj}})$ and $\tilde{X}'_{njk} = (1, \tilde{X}'_{n1jk}, \ldots, \tilde{X}'_{npjk})$, $X'_{njk} = (1, X'_{n1jk}, \ldots, X'_{npjk})$ for $k = 1, \ldots, L_{nj}$ contain $p$ components of the $k$th element in bin $j$, and

$$\tilde{X}'_{nj} = (\tilde{X}'^T_{nj1}, \ldots, \tilde{X}'^T_{njL_{nj}})^T_{L_{nj} \times (p+1)}, \qquad X'_{nj} = (X'^T_{nj1}, \ldots, X'^T_{njL_{nj}})^T_{L_{nj} \times (p+1)}.$$

Then we can express the least squares estimates of the multiple regression of the observable data falling in the $j$th bin $B_{nj}$ as

$$(9) \qquad \hat{\beta}_{nj}^T = (\hat{\beta}_{n0j}, \ldots, \hat{\beta}_{npj})^T = (\tilde{X}'^T_{nj} \tilde{X}'_{nj})^{-1} \tilde{X}'^T_{nj} \tilde{Y}'^T_{nj},$$

leading to the parameter estimates $\hat{\gamma}_{n0}$ and $\hat{\gamma}_{nr}$ given in (6) and (7), respectively, where $\bar{\tilde{X}}_{nr} = n^{-1} \sum_{i=1}^{n} \tilde{X}_{nri}$ and $\bar{\tilde{X}}'_{nrj} = L_{nj}^{-1} \sum_{k=1}^{L_{nj}} \tilde{X}'_{nrjk}$.

Let $\tilde{\gamma}_{nj}$ be the least squares estimates of the multiple regression of the unobservable data falling into $B_{nj}$, that is,

$$(10) \qquad \tilde{\gamma}_{nj}^T = (\tilde{\gamma}_{n0j}, \ldots, \tilde{\gamma}_{npj})^T = (X'^T_{nj} X'_{nj})^{-1} X'^T_{nj} Y'^T_{nj}.$$

This quantity is not estimable, but will be used in the proof of the main results.

For the estimates given in (6) and (7) to be well defined, the least squares estimate $\hat{\beta}_{nj}$ must exist for each bin $B_{nj}$. This requires that the inverse of $\tilde{X}'^T_{nj} \tilde{X}'_{nj}$ is well defined, that is, $\det(\tilde{X}'^T_{nj} \tilde{X}'_{nj}) \neq 0$. Correspondingly, $\tilde{\gamma}_{nj}$ will exist under the condition that $\det(X'^T_{nj} X'_{nj}) \neq 0$. Define the events

$$(11) \quad \begin{aligned} \tilde{A}_n &= \left\{\omega \in \Omega : \inf_j |\det(L_{nj}^{-1} \tilde{X}'^T_{nj} \tilde{X}'_{nj})| > \zeta \text{ and } \min_j L_{nj} > p\right\}, \\ A_n &= \left\{\omega \in \Omega : \inf_j |\det(L_{nj}^{-1} X'^T_{nj} X'_{nj})| > \zeta \text{ and } \min_j L_{nj} > p\right\}, \end{aligned}$$

where $\zeta = \min\{\rho/2, [\inf_j(\phi_1^2(U'^*_{nj}), \ldots, \phi_p^2(U'^*_{nj}))]^p \rho/2\}$, $\rho$ is as defined in (C5), $U'^*_{nj} = L_{nj}^{-1} \sum_{k=1}^{L_{nj}} U'_{njk}$ is the average of the $U$'s in $B_{nj}$ and $(\Omega, \mathcal{F}, P)$ is the underlying probability space. On event $\tilde{A}_n$, $\hat{\gamma}_{n0}$ and $\hat{\gamma}_{nr}$ given in (6) and (7), and on event $A_n$, $\tilde{\gamma}_{nj}$ given in (10) are well defined, respectively. Event $E_n$ in Theorems 1 and 2 is defined to be the intersection of $A_n$ and $\tilde{A}_n$, that is,

$$E_n = A \cap \tilde{A}_n. \tag{12}$$

It is shown in Appendix A.3 that $P(E_n) \to 1$ as $n \to \infty$.

THEOREM 2. *Under the technical conditions* (C1)–(C7) *in Section* 6, *on event* $E_n$ [*defined in* (12)] *with* $P(E_n) \to 1$ *as* $n \to \infty$,

$$\hat{\sigma}_{nr}^2 \xrightarrow{p} \sigma_r^2, \qquad 0 \le r \le p,$$

*where*

$$\hat{\sigma}_{n0}^2 = \left(\sum_{j=1}^m \frac{L_{nj}}{n}\hat{\beta}_{n0j}^2 - \hat{\gamma}_{n0}^2\right)$$

$$+ \left\{\frac{1}{n}\sum_{j=1}^m \sum_{k=1}^{L_{nj}} (\tilde{Y}'_{njk} - \hat{\beta}_{n0j} - \hat{\beta}_{n1j}\tilde{X}'_{n1jk} - \cdots - \hat{\beta}_{npj}\tilde{X}'_{npjk})^2\right\}$$

$$\times \left\{\sum_{j=1}^m \frac{L_{nj}}{n}(L_{nj}^{-1}\tilde{X}'^T_{nj}\tilde{X}'_{nj})_{11}^{-1}\right\},$$

$$\hat{\sigma}_{nr}^2 = \left[\frac{1}{n}\sum_{j=1}^m \hat{\beta}_{nrj}^2 \sum_{k=1}^{L_{nj}} \tilde{X}'^2_{nrjk} + \hat{\gamma}_{nr}^2 \bar{\tilde{X}}_{nr}^2 - 2\frac{\hat{\gamma}_{nr}}{n}\sum_{j=1}^m \hat{\beta}_{nrj}\sum_{k=1}^{L_{nj}}\tilde{X}'^2_{nrjk} + \hat{\gamma}_{nr}^2 s_{\tilde{X}_r}^2\right.$$

$$+ \left\{\frac{1}{n}\sum_{j=1}^m \sum_{k=1}^{L_{nj}}(\tilde{Y}'_{njk} - \hat{\beta}_{n0j} - \hat{\beta}_{n1j}\tilde{X}'_{n1jk} - \cdots - \hat{\beta}_{npj}\tilde{X}'_{npjk})^2\right\}$$

$$\left. \times \left\{\sum_{j=1}^m \frac{L_{nj}}{n}\bar{\tilde{X}}_{nrj}^2(L_{nj}^{-1}\tilde{X}'^T_{nj}\tilde{X}'_{nj})_{rr}^{-1}\right\}\right]\Big/ \bar{\tilde{X}}_{nr}^2,$$

$$1 \le r \le p,$$

*and* $s_{\tilde{X}_r}^2 = (n-1)^{-1}\sum_{i=1}^n (\tilde{X}_{nri} - \bar{\tilde{X}}_{nr})^2$.

REMARK. These proposed variance estimates are motivated by the identifiability conditions, the definition of the smooth varying coefficient functions given in (5), Lemma A.3 and Lemma A.4(a). Using the consistency of $\hat{\beta}_{nrj}$ for the value of the function $\beta_r$ at the midpoint of the $j$th bin and the definitions of $\tilde{Y}'_{njk}$ and $\tilde{X}'_{nrjk}$, we target the quantities $\sigma^2 E\{\psi(U)\}$, $\gamma_0^2 E\{\psi^2(U)\}$, $\gamma_r^2 E(X_r^2)E\{\psi^2(U)\}$ and $\gamma_r^2 E\{\phi_r(U)\psi(U)\}E(X_r^2)$ with the estimators $n^{-1} \times \sum_{j=1}^m \sum_{k=1}^{L_{nj}}(\tilde{Y}'_{njk} - \hat{\beta}_{n0j} - \hat{\beta}_{n1j}\tilde{X}'_{n1jk} - \cdots - \hat{\beta}_{npj}\tilde{X}'_{npjk})^2$, $\sum_{j=1}^m n^{-1}L_{nj}\hat{\beta}_{n0j}^2$, $n^{-1}\sum_{j=1}^m \hat{\beta}_{nrj}^2\sum_{k=1}^{L_{nj}}\tilde{X}'^2_{nrjk}$ and $n^{-1}\hat{\gamma}_{nr}\sum_{j=1}^m \hat{\beta}_{nrj}\sum_{k=1}^{L_{nj}}\tilde{X}'^2_{nrjk}$,

respectively. Furthermore, relying mainly on Lemmas A.3 and A.4(a), we target $(\mathcal{X}^{-1})_{11}$ and $\{E(X_r)\}^2(\mathcal{X}^{-1})_{rr}$ with $\sum_{j=1}^{m} n^{-1} L_{nj} (L_{nj}^{-1} \tilde{X}_{nj}'^T \tilde{X}_{nj}')_{11}^{-1}$ and $\sum_{j=1}^{m} n^{-1} L_{nj} \bar{\tilde{X}}_{nrj}^2 \times (L_{nj}^{-1} \tilde{X}_{nj}'^T \tilde{X}_{nj}')_{rr}^{-1}$, respectively.

**5. Applications and Monte Carlo study.** Under the technical conditions (C1)–(C7) in Section 6,

$$
(13) \qquad \frac{\sqrt{n}}{\sigma_r}(\hat{\gamma}_{nr} - \gamma_r) \xrightarrow{\mathcal{D}} \mathbb{N}(0,1), \qquad 0 \leq r \leq p \text{ as } n \to \infty.
$$

Using the consistent estimate $\hat{\sigma}_{nr}^2$ of $\sigma_r^2$ proposed in Theorem 2, it follows from (13) and Slutsky's theorem that

$$
\frac{\sqrt{n}}{\hat{\sigma}_{nr}}(\hat{\gamma}_{nr} - \gamma_r) \xrightarrow{\mathcal{D}} \mathbb{N}(0,1), \qquad 0 \leq r \leq p,
$$

so that an approximate $(1-\alpha)$ asymptotic confidence interval for $\gamma_r$ has the endpoints

$$
(14) \qquad \hat{\gamma}_{nr} \pm z_{\alpha/2} \frac{\hat{\sigma}_{nr}}{\sqrt{n}}.
$$

Here $z_{\alpha/2}$ is the $(1-\alpha/2)$th quantile of the standard Gaussian distribution.

5.1. *Application to creatinine data.* An observational study in which various laboratory and patient data were analyzed for patients with end-stage renal disease is described in [7]. To illustrate our methods, we analyzed a similar but much smaller set of data and note that our analysis does not provide inference for the data in [7]. Variables include serum creatinine level ($CRT$), cholesterol level ($CH$), serum albumin level ($ALB$) and body mass index ($BMI$), measured for $n = 508$ subjects. Creatinine is a protein produced by muscle and released into the blood. Since the amount produced is relatively stable, the creatinine level in the serum is determined by the rate at which it is removed, and is therefore an important indicator of renal function. We analyze the dependence of serum creatinine (response) on cholesterol level and serum albumin (predictors). An unadjusted approach would be to fit the multiple regression model $CRT = \gamma_0 + \gamma_1 CH + \gamma_2 ALB + e$, where $e$ is an error term, usually by least squares. Body mass index ($BMI$) is defined as weight/height$^2$ and is known to affect both the response and the predictors. This provides the motivation to adjust for this influence by means of the CAR model (4), (5), using body mass index as the confounder $U$.

The parameters $\gamma_0$, $\gamma_1$ and $\gamma_2$ were estimated by the CAR algorithm and the results were compared to the estimates obtained from the least squares regression of the observed $CRT$ on observed $CH$ and $ALB$. The estimates and the approximate 95% asymptotic confidence intervals for the regression parameters obtained through both methods are displayed in Table 1.

TABLE 1
*Parameter estimates for the regression model $CRT = \gamma_0 + \gamma_1 CH + \gamma_2 ALB + e$, obtained by least squares regression of $\tilde{Y} = CRT$ (serum creatinine level) on $\tilde{X}_1 = CH$ (cholesterol level) and $\tilde{X}_2 = ALB$ (serum albumin level), and alternatively by covariate adjusted regression, for $n = 508$ subjects*

| | Least sq. reg. | | | Covariate adj. reg. | | |
|---|---|---|---|---|---|---|
| Coefficients | Lower b. | Estimate | Upper b. | Lower b. | Estimate | Upper b. |
| Intercept | 1.2715 | 4.3685 | 7.4656 | 0.3679 | 3.9987 | 7.6296 |
| CH | −0.0106 | −0.0041 | 0.0023 | −0.0154 | −0.0082 | −0.0009 |
| ALB | 1.1819 | 1.9729 | 2.7639 | 1.3065 | 2.2532 | 3.2000 |

Confidence intervals at the 95% level were obtained by the standard $t$-statistic for least squares regression and by the proposed asymptotic intervals (14) for CAR, respectively.

The approximate confidence intervals for CAR estimates were obtained as proposed in (14). The scatter-plots of the raw estimates $(\hat{\beta}_{nr1}, \ldots, \hat{\beta}_{nrm})$ (9) versus midpoints of the bins $(B_{n1}, \ldots, B_{nm})$ are shown in Figure 1 for $r = 0, 1, 2$.

The implementation of the binning algorithm allows for merging of sparsely populated bins. Bin widths were chosen such that there are at least $(p+1)$ points, enough to fit the linear regression with $(p-1)$ predictors in each bin. Where there were bins with less than $(p+1)$ elements, such bins were merged with neighboring bins. For this example with $n = 508$, the average number of points per bin was 14, yielding a total of 34 bins after merging.

For least squares regression, *CH* was not found significant at the usual 5% level, while *ALB* was found to be significant. When applying the CAR method, *CH* and *ALB* were both significant. As *BMI* increases, the slope parameter of serum albumin level increases exponentially, while the negative slope parameter of cholesterol level declines slightly. Adjusting for different *BMI* levels across patients, both serum albumin level and cholesterol level seem to play a significant role for the serum creatinine level. The effects of *BMI* are thus masking the true overall negative effect that *CH* has on *CRT* in the unadjusted regression equation.

5.2. *Monte Carlo simulation.* The confounding covariate $U$ was simulated from Uniform$(2,6)$. The underlying unobserved multiple regression model was

$$(15) \qquad Y = 4 - X_1 + 0.3X_2 + 3X_3 + e,$$

where $X_1 \sim \mathcal{N}(1.5, 0.7)$, $X_2 \sim \mathcal{N}(1, 1.2)$, $X_3 \sim \mathcal{N}(0.5, 1)$ and $e \sim \mathcal{N}(0, 0.3)$. The distortion functions were chosen as $\psi(U) = (U+3)/7$, $\phi_1(U) = (U+1)^2/26.3333$, $\phi_2(U) = (U+10)/14$ and $\phi_3(U) = (U+2)^2/37.3333$, satisfying the identifiability conditions. We conducted 1000 Monte Carlo runs with sample sizes 100, 400 and 1600. For each run approximate 95% asymptotic confidence intervals were formed for the regression parameters by plugging in the estimates $\hat{\sigma}^2_{nr}$, $r = 0, \ldots, p$, given in Theorem 2, into (14). The estimated coverage fractions and mean interval lengths for these confidence intervals are given in Table 2. The estimated noncoverage fractions are seen to get very close to the target value 0.05 as sample size increases, and the estimated interval lengths are sharply decreasing.

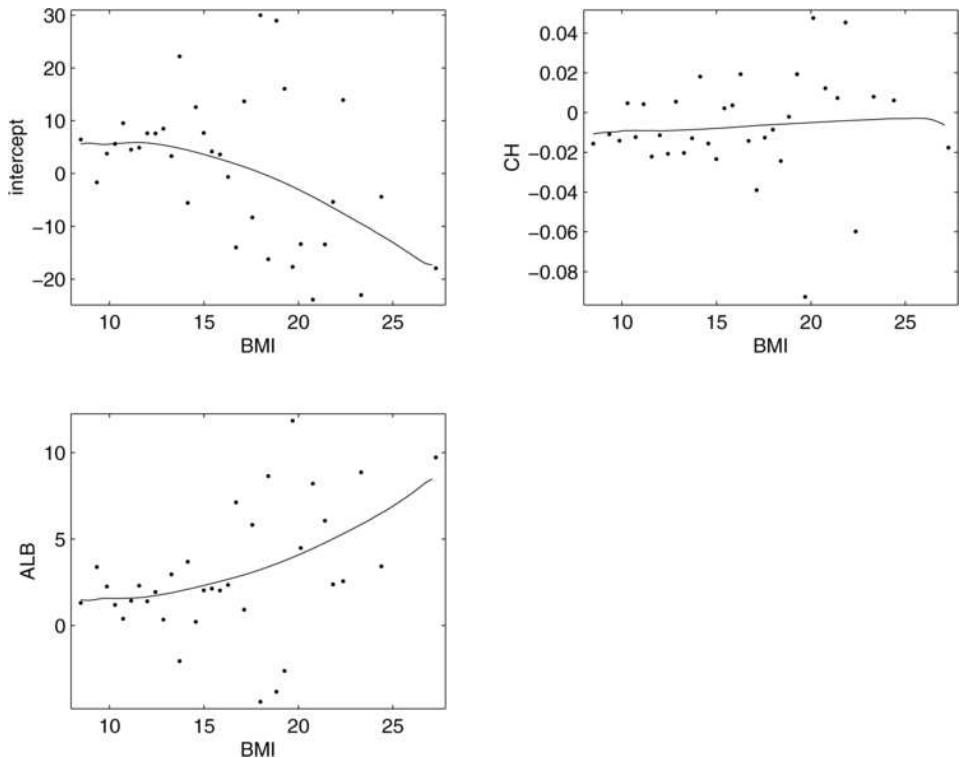

FIG. 1. *Scatter-plots of the raw estimates $(\hat{\beta}_{nr1}, \ldots, \hat{\beta}_{nrm})$ versus midpoints of the bins $(B_{n1}, \ldots, B_{nm})$ for $r = 0$ (top left panel) and $r = 1$ (top right panel) and $r = 2$ (bottom left panel) in the CAR model $CRT = \beta_0(BMI) + \beta_1(BMI)CH + \beta_2(BMI)ALB + \varepsilon(BMI)$. Local polynomial smooth curves have been fitted through the scatter-plots using cross-validation bandwidth choices of $h = 8, 7, 7$, respectively, for $r = 0, 1, 2$. CRT = serum creatinine level, CH = cholesterol level, ALB = serum albumin level and BMI = body mass index. Sample size is* 508, *and the number of bins formed is* 34.

We have also carried out simulations to study the effects of different choices of $m$, the total number of bins, on the mean square error of the CAR estimates. Under the rate conditions on $m$ given in Section 3, the estimates are found to be sufficiently robust regarding different choices of $m$.

**6. Proofs of the main results.** While the main steps in the proofs of the two theorems are given here, the auxiliary results for these proofs are deferred to the Appendix, where they are listed as Lemmas A.1–A.4. We introduce some technical conditions:

(C1) The covariate $U$ is bounded below and above, $-\infty < a \leq U \leq b < \infty$ for real numbers $a < b$. The density $f(u)$ of $U$ satisfies $\inf_{a \leq u \leq b} f(u) > c_1 > 0$, $\sup_{a \leq u \leq b} f(u) < c_2 < \infty$ for real $c_1, c_2$, and is uniformly Lipschitz continuous, that is, there exists a real number $M$ such that $\sup_{a \leq u \leq b} |f(u+c) - f(u)| \leq M|c|$ for any real number $c$.
(C2) The variables $(e, U, X_r)$ are mutually independent for $r = 1, \ldots, p$.
(C3) For the predictors, $\sup_{1 \leq i \leq n, 1 \leq r \leq p} |X_{nri}| \leq B$ for some bound $B \in \mathbb{R}$.
(C4) Contamination functions $\psi(\cdot)$ and $\phi_r(\cdot)$, $1 \leq r \leq p$, are twice continuously differentiable, satisfying

$$E\psi(U) = 1, \qquad E\phi_r(U) = 1, \qquad \phi_r(\cdot) > 0, 1 \leq r \leq p.$$

(C5) As $n \to \infty$, $\frac{1}{n} X^T X \xrightarrow{p} \mathcal{X}$, where $\mathcal{X}$, the limiting $(p+1) \times (p+1)$-matrix, is nonsingular, that is, $\rho = |\det(\mathcal{X})| > 0$.

These are mild conditions that are satisfied in most practical situations. Bounded covariates are standard in asymptotic theory for least squares regression, as are conditions (C2) and (C5) (see [8]). The identifiability conditions stated in (C4) are equivalent to

$$E(\tilde{Y}|X) = E(Y|X), \qquad E(\tilde{X}_r|X_r) = X_r.$$

This means that the confounding of $Y$ by $U$ does not change the mean regression function. Some further technical conditions will be introduced in Appendix A.1; these are required to prove the auxiliary lemmas in the Appendix.

For two matrices of the same dimension, let $A \boxdot B$ denote the Hadamard product, where $A \boxdot B$ is also of the same dimension with $(i, j)$th element

TABLE 2
*Coverage (in percent) and mean interval length for the approximate 95% asymptotic confidence intervals formed for the parameters of the regression model* (15)

|   | $\gamma_0$ | | $\gamma_1$ | | $\gamma_2$ | | $\gamma_3$ | |
|---|---|---|---|---|---|---|---|---|
| $n$ | Coverage | Length | Coverage | Length | Coverage | Length | Coverage | Length |
| 100  | 90.7 | 0.56 | 90.4 | 0.32 | 91.7 | 0.20 | 96.6 | 0.73 |
| 400  | 93.4 | 0.21 | 94.1 | 0.11 | 93.4 | 0.06 | 95.5 | 0.30 |
| 1600 | 94.2 | 0.10 | 95.2 | 0.05 | 94.7 | 0.03 | 95.0 | 0.14 |

The values were obtained from 1000 Monte Carlo runs. The average number of points per bin was 5, 16 and 32 for sample sizes 100, 400 and 1600.

equal to the product of the $(i,j)$th elements of $A$ and $B$.

PROOF OF THEOREM 1. By Lemma A.4(b) and properties (b), (c), (e), (f) given in Appendix A.3, it holds that

$$(16) \quad \sup_j |(L_{nj}^{-1}\tilde{X}_{nj}^{\prime T}\tilde{Y}_{nj}^{\prime}) - \{\Delta \boxdot (L_{nj}^{-1}X_{nj}^{\prime T}Y_{nj}^{\prime})\}| = O_p(m^{-1})\mathbf{1}_{(p+1)\times 1},$$

where $(L_{nj}^{-1}\tilde{X}_{nj}^{\prime T}\tilde{Y}_{nj}^{\prime}) = (L_{nj}^{-1}\sum_k \tilde{Y}_{njk}^{\prime}, L_{nj}^{-1}\sum_k \tilde{Y}_{njk}^{\prime}\tilde{X}_{n1jk}^{\prime},\ldots,L_{nj}^{-1}\sum_k \tilde{Y}_{njk}^{\prime}\tilde{X}_{npjk}^{\prime})^T$, $(L_{nj}^{-1}X_{nj}^{\prime T}Y_{nj}^{\prime}) = (L_{nj}^{-1}\sum_k Y_{njk}^{\prime}, L_{nj}^{-1}\sum_k Y_{njk}^{\prime}X_{n1jk}^{\prime},\ldots,L_{nj}^{-1}\sum_k Y_{njk}^{\prime}X_{npjk}^{\prime})^T$, $\Delta = \{\psi(U_{nj}^{\prime *}),\psi(U_{nj}^{\prime *})\phi_1(U_{nj}^{\prime *}),\ldots,\psi(U_{nj}^{\prime *})\phi_p(U_{nj}^{\prime *})\}^T$ and $\mathbf{1}_{(p+1)\times 1}$ denotes a $(p+1)\times 1$ vector of 1's. Under event $E_n$, Lemma A.3 and (16) imply that

$$(17) \quad \sup_j \begin{vmatrix} \hat{\beta}_{n0j} & - & \psi(U_{nj}^{\prime *})\tilde{\gamma}_{n0j} \\ \hat{\beta}_{n1j} & - & \{\psi(U_{nj}^{\prime *})/\phi_1(U_{nj}^{\prime *})\}\tilde{\gamma}_{n1j} \\ \vdots & & \vdots \\ \hat{\beta}_{npj} & - & \{\psi(U_{nj}^{\prime *})/\phi_p(U_{nj}^{\prime *})\}\tilde{\gamma}_{npj} \end{vmatrix} = O_p(m^{-1})\mathbf{1}_{(p+1)\times 1},$$

where $\tilde{\gamma}_{nj}$ is as defined in (10). First consider the case $r=0$. Using (17),

$$\sqrt{n}(\hat{\gamma}_{n0} - \gamma_0)$$
$$= \sqrt{n}\left(\sum_j^m \frac{L_{nj}}{n}\hat{\beta}_{n0j} - \gamma_0\right)$$
$$= \sum_j^m \frac{L_{nj}}{\sqrt{n}}\psi(U_{nj}^{\prime *})\tilde{\gamma}_{n0j} - \sqrt{n}\gamma_0 + O_p\left(\frac{\sqrt{n}}{m}\right)$$
$$= \sum_j^m \frac{L_{nj}}{\sqrt{n}}\psi(U_{nj}^{\prime *})[\gamma_0 + \{(X_{nj}^{\prime T}X_{nj}^{\prime})^{-1}X_{nj}^{\prime T}e_{nj}^{\prime}\}_1] - \sqrt{n}\gamma_0 + O_p\left(\frac{\sqrt{n}}{m}\right).$$

By property (b), Lemma A.4(a), (b) and substituting $L_{nj}^{-1}\sum_k\{(L_{nj}^{-1}X_{nj}^{\prime T}X_{nj}^{\prime})^{-1}\times X_{nj}^{\prime T}\}_{1k}e_{njk}^{\prime}$ for $\{(X_{nj}^{\prime T}X_{nj}^{\prime})^{-1}X_{nj}^{\prime T}e_{nj}^{\prime}\}_1$, $\sqrt{n}(\hat{\gamma}_{n0} - \gamma_0)$ further simplifies to

$$(18) \quad \sum_{j=1}^m \sum_{k=1}^{L_{nj}} \left[\frac{\gamma_0\psi(U_{njk}^{\prime})}{\sqrt{n}} + \frac{\psi(U_{njk}^{\prime})e_{njk}^{\prime}}{\sqrt{n}}\{(L_{nj}^{-1}X_{nj}^{\prime T}X_{nj}^{\prime})^{-1}X_{nj}^{\prime T}\}_{1k}\right]$$
$$- \sqrt{n}\gamma_0 + O_p\left(\frac{\sqrt{n}}{m}\right).$$

Since the above sum is over all bins indexed by $j$, and over all points within the bins indexed by $k$, it is equal to the sum over all data points indexed by $i$, summed up in a random order. We introduce notation where $X'_{nj(i)}$ refers to the matrix $X'_{nj}$ and $L_{nj(i)}$ refers to the number of points in the $j$th bin such that $U_{ni} \in B_{nj}$, and $\{(L^{-1}_{nj(i)} X'^T_{nj(i)} X'_{nj(i)})^{-1} X'^T_{nj(i)}\}_{rk(i)}$ is the $(r,k)$th element of the matrix $\{(L^{-1}_{nj} X'^T_{nj} X'_{nj})^{-1} X'^T_{nj}\}$ for $1 \leq r \leq p$, where $U_{ni} = U'_{njk}$ is the $k$th element in the ordered sample $(U'_{nj1}, \ldots, U'_{njL_{nj}}) \in B_{nj}$. Thus (18) is equal to

$$(19) \quad \sum_{\substack{i=1 \\ j,k}}^{n} \left[ \frac{\gamma_0 \psi(U_{ni})}{\sqrt{n}} + \frac{\psi(U_{ni}) e_{ni}}{\sqrt{n}} \{(L^{-1}_{nj(i)} X'^T_{nj(i)} X'_{nj(i)})^{-1} X'^T_{nj(i)}\}_{1k(i)} - \frac{\gamma_0}{\sqrt{n}} \right]$$
$$+ O_p\left(\frac{\sqrt{n}}{m}\right).$$

The term $\sqrt{n}(\hat{\gamma}_{n0} - \gamma_0)$ is asymptotically equivalent to

$$S_{n0t} = \sum_{\substack{i=1 \\ j,k}}^{t} \left[ \frac{\gamma_0 \psi(U_{ni})}{\sqrt{n}} + \frac{\psi(U_{ni}) e_{ni}}{\sqrt{n}} \{(L^{-1}_{nj(i)} X'^T_{nj(i)} X'_{nj(i)})^{-1} X'^T_{nj(i)}\}_{1k(i)} - \frac{\gamma_0}{\sqrt{n}} \right]$$
$$= \sum_{i=1}^{t} Z_{n0i},$$

since $m/\sqrt{n} \to \infty$ as $n \to \infty$ makes the term $O_p(\sqrt{n}/m)$ negligible.

Let $F_{n0t}$ be the $\sigma$-field generated by $\{e_{n1}, \ldots, e_{nt}, U_{n1}, \ldots, U_{nt}, L_{nj(1)}, \ldots, L_{nj(t)}, X'_{nj(1)}, \ldots, X'_{nj(t)}\}$. Then $\{S_{n0t} = \sum_{i=1}^{t} Z_{n0i}, F_{n0t}, 1 \leq t \leq n\}$ is a mean-zero martingale for $n \geq 1$, since $E(S_{n0t}) = 0$, $E(S_{n0,t+1}|F_{n0t}) = S_{n0t}$ and $S_{n0t}$ is adapted to $F_{n0t}$. Since the $\sigma$-fields are nested, that is, $F_{n0t} \subseteq F_{n0,t+1}$ for all $t \leq n$, using Lemma A.1, $S_{n0n} \to \mathbb{N}(0, \sigma_0^2)$ in distribution ([9], Theorem 2.3 and subsequent discussion), and Theorem 1 follows for $r = 0$.

Next we show

$$(20) \quad \sqrt{n} \begin{pmatrix} \sum_{j=1}^{m} \frac{L_{nj}}{n} \hat{\beta}_{nrj} \bar{\tilde{X}}'_{nrj} - \gamma_r E(X_r) \\ \sum_{j=1}^{m} \frac{L_{nj}}{n} \bar{\tilde{X}}'_{nrj} - E(X_r) \end{pmatrix} \xrightarrow{\mathcal{D}} \mathbb{N}_2(\underline{0}, \Sigma_r).$$

The asymptotic normality of $\sqrt{n}(\hat{\gamma}_{nr} - \gamma_r)$ for $r = 1, \ldots, p$ will follow from this with a simple application of the $\delta$-method, since $\hat{\gamma}_{nr} = (\sum_{j=1}^{m} L_{nj} n^{-1} \hat{\beta}_{nrj} \bar{\tilde{X}}'_{nrj}) / (\sum_{j=1}^{m} L_{nj} n^{-1} \bar{\tilde{X}}'_{nrj})$ as defined in (7). By the Cramér–Wald device it is enough to show the asymptotic normality of $\sqrt{n}[a\{\sum_{j=1}^{m} L_{nj} n^{-1} \hat{\beta}_{nrj} \bar{\tilde{X}}'_{nrj} - \gamma_r E(X_r)\} + b\{\sum_{j=1}^{m} L_{nj} n^{-1} \bar{\tilde{X}}'_{nrj} - E(X_r)\}]$ for real $a, b$, and (20) will follow.

Using (17), properties (b), (c), Lemma A.4(a), (b) and substituting $L_{nj}^{-1} \times \sum_k \{(L_{nj}^{-1} X_{nj}'^T X_{nj}')^{-1} X_{nj}'^T\}_{rk} e_{njk}'$ for $\{(X_{nj}'^T X_{nj}')^{-1} X_{nj}'^T e_{nj}'\}_r$, we have

$$\sum_{j=1}^m \frac{L_{nj}}{n} \hat{\beta}_{nrj} \bar{\bar{X}}_{nrj}' = \sum_{j=1}^m \frac{L_{nj}}{n} \psi(U_{nj}'^*) \bar{X}_{nrj}' [\gamma_r + \{(X_{nj}'^T X_{nj}')^{-1} X_{nj}'^T e_{nj}'\}_r] + O_p(m^{-1})$$

$$= \sum_{j=1}^m \sum_{k=1}^{L_{nj}} \left[ \frac{\gamma_r}{n} \psi(U_{njk}') X_{nrjk}' \right.$$

$$\left. + \frac{\bar{X}_{nrj}'}{n} \psi(U_{njk}') e_{njk}' \{(L_{nj}^{-1} X_{nj}'^T X_{nj}')^{-1} X_{nj}'^T\}_{rk} \right]$$

$$+ O_p(m^{-1})$$

and

$$\sum_{j=1}^m \frac{L_{nj}}{n} \bar{\bar{X}}_{nrj}' = \sum_{j=1}^m \sum_{k=1}^{L_{nj}} \frac{1}{n} \phi_r(U_{njk}') X_{nrjk}' + O_p(m^{-1}).$$

Thus using the same notation as in (19), it holds that

$$\sqrt{n} \left[ a \left\{ \sum_{j=1}^m \frac{L_{nj}}{n} \hat{\beta}_{nrj} \bar{\bar{X}}_{nrj}' - \gamma_r E(X_r) \right\} + b \left\{ \sum_{j=1}^m \frac{L_{nj}}{n} \bar{\bar{X}}_{nrj}' - E(X_r) \right\} \right]$$

$$= \sum_{\substack{i=1 \\ j,k}}^n \left[ a \frac{\gamma_r}{\sqrt{n}} \psi(U_{ni}) X_{nri} \right.$$

$$+ a \frac{\bar{X}_{nrj(i)}'}{\sqrt{n}} \psi(U_{ni}) e_{ni} \{(L_{nj(i)}^{-1} X_{nj(i)}'^T X_{nj(i)}')^{-1} X_{nj(i)}'^T\}_{rk(i)}$$

$$\left. - a \frac{\gamma_r}{\sqrt{n}} E(X_r) + \frac{b}{\sqrt{n}} \phi_r(U_{ni}) X_{nri} - b \frac{E(X_r)}{\sqrt{n}} \right]$$

$$+ O_p \left( \frac{\sqrt{n}}{m} \right),$$

where $\bar{X}_{nrj(i)}' = L_{nj}^{-1} \sum_{k=1}^{L_{nj(i)}} X_{nrj(i)k}'$. Since $O_p(\sqrt{n}/m)$ is asymptotically negligible, the above term is asymptotically equivalent to

$$S_{nrt} = \sum_{\substack{i=1 \\ j,k}}^n \left[ a \frac{\gamma_r}{\sqrt{n}} \psi(U_{ni}) X_{nri} \right.$$

$$+ a \frac{\bar{X}_{nrj(i)}'}{\sqrt{n}} \psi(U_{ni}) e_{ni} \{(L_{nj(i)}^{-1} X_{nj(i)}'^T X_{nj(i)}')^{-1} X_{nj(i)}'^T\}_{rk(i)}$$

$$- a\frac{\gamma_r}{\sqrt{n}}E(X_r) + \frac{b}{\sqrt{n}}\phi_r(U_{ni})X_{nri} - b\frac{E(X_r)}{\sqrt{n}}\Bigg]$$

$$= \sum_{i=1}^{t} Z_{nri}.$$

Let $F_{nrt}$ be the $\sigma$-field generated by $\{e_{n1},\ldots,e_{nt},U_{n1},\ldots,U_{nt},L_{nj(1)},\ldots,L_{nj(t)},X'_{nj(1)},\ldots,X'_{nj(t)}\}$. Then it is easy to check that $\{S_{nrt} = \sum_{i=1}^{t} Z_{nri}, F_{nrt}, 1 \leq t \leq n\}$ is a mean-zero martingale for $n \geq 1$. Since the $\sigma$-fields are nested, that is, $F_{nrt} \subseteq F_{nr,t+1}$ for all $t \leq n$, using Lemma A.2, $S_{nrn} \xrightarrow{\mathcal{D}} \mathbb{N}(0,(a,b)\Sigma_r(a,b)^T)$. Thus, it also follows by a simple application of the $\delta$-method that $\sqrt{n}(\hat{\gamma}_{nr} - \gamma_r) \xrightarrow{p} \mathbb{N}(0,\sigma_r^2)$ for $r = 1,\ldots,p$, where $\sigma_r^2$ is as defined in Theorem 1. $\square$

PROOF OF THEOREM 2. Using Lemma A.4(a) and (b), it holds on event $A_n$ that

$$\sup_{j} |\tilde{\gamma}_{nj} - \underline{\gamma}| = o_p(1)\mathbf{1}_{(p+1)\times 1}, \tag{21}$$

where $\underline{\gamma} = (\gamma_0, \gamma_1, \ldots, \gamma_p)^T$. Using (21) and (17),

$$\sup_{j} \begin{vmatrix} \hat{\beta}_{n0j} & - & \psi(U'^*_{nj})\gamma_0 \\ \hat{\beta}_{n1j} & - & \{\psi(U'^*_{nj})/\phi_1(U'^*_{nj})\}\gamma_1 \\ \vdots & & \vdots \\ \hat{\beta}_{npj} & - & \{\psi(U'^*_{nj})/\phi_p(U'^*_{nj})\}\gamma_p \end{vmatrix} = o_p(1)\mathbf{1}_{(p+1)\times 1}. \tag{22}$$

By (22), properties (b), (c), (d), boundedness considerations and the law of large numbers,

$$\sum_{j=1}^{m} \frac{L_{nj}}{n}\hat{\beta}_{n0j}^2 = \gamma_0^2 \sum_{j=1}^{m} \frac{L_{nj}}{n}\{\psi^2(U'^*_{nj}) + o_p(1)\}$$

$$= \frac{\gamma_0^2}{n}\sum_{i=1}^{n}\psi^2(U_{ni}) + o_p(1) = \gamma_0^2 E\{\psi^2(U)\} + o_p(1),$$

$$\frac{1}{n}\sum_{j=1}^{m}\sum_{k=1}^{L_{nj}}(\tilde{Y}'_{njk} - \hat{\beta}_{n0j} - \hat{\beta}_{n1j}\tilde{X}'_{n1jk} - \cdots - \hat{\beta}_{npj}\tilde{X}'_{npjk})^2$$

$$= \frac{1}{n}\sum_{j=1}^{m}\sum_{k=1}^{L_{nj}}\bigg\{\psi(U'^*_{nj})e'_{njk} + \delta_{n0jk}Y'_{njk} - \gamma_1\frac{\psi(U'^*_{nj})}{\phi_1(U'^*_{nj})}\delta_{n1jk}X'_{n1jk} - \cdots$$

$$- \gamma_p\frac{\psi(U'^*_{nj})}{\phi_p(U'^*_{nj})}\delta_{npjk}X'_{npjk} + o_p(1)\bigg\}^2$$

$$= \frac{1}{n} \sum_{i=1}^{n} \psi(U_{ni}^2) e_{ni}^2 + o_p(1) = \sigma^2 E\{\psi^2(U)\} + o_p(1),$$

$$\frac{1}{n} \sum_{j=1}^{m} \hat{\beta}_{nrj}^2 \sum_{k=1}^{L_{nj}} \tilde{X}_{nrjk}'^2 = \frac{\gamma_r^2}{n} \sum_{j=1}^{m} \psi^2(U_{nj}'^*) \sum_{k=1}^{L_{nj}} X_{nrjk}'^2 + o_p(1)$$

$$= \frac{\gamma_r^2}{n} \sum_{i=1}^{n} \psi^2(U_{ni}) X_{nri}^2 + o_p(1) = \gamma_r^2 E\{\psi^2(U)\} E(X_r^2) + o_p(1)$$

and

$$\frac{1}{n} \sum_{j=1}^{m} \hat{\beta}_{nrj} \sum_{k=1}^{L_{nj}} \tilde{X}_{nrjk}'^2 = \frac{\gamma_r}{n} \sum_{j=1}^{m} \psi(U_{nj}'^*) \phi_r(U_{nj}'^*) \sum_{k=1}^{L_{nj}} X_{nrjk}'^2 + o_p(1)$$

$$= \frac{\gamma_r}{n} \sum_{i=1}^{n} \psi(U_{ni}) \phi_r(U_{ni}) X_{nri}^2 + o_p(1)$$

$$= \gamma_r E\{\psi(U) \phi_r(U)\} E(X_r^2) + o_p(1),$$

where $\delta_{n0jk}$ and $\delta_{nrjk}$ are as defined in Appendix A.3. Using Lemma A.3, Lemma A.4(a) and (31),

$$\sum_{j=1}^{m} \frac{L_{nj}}{n} \left( \frac{1}{L_{nj}} \tilde{X}_{nj}' \tilde{X}_{nj}' \right)^{-1}_{11} \xrightarrow{p} (\mathcal{X}^{-1})_{11},$$

$$\sum_{j=1}^{m} \frac{L_{nj}}{n} \bar{\tilde{X}}_{nrj}^2 \left( \frac{1}{L_{nj}} \tilde{X}_{nj}' \tilde{X}_{nj}' \right)^{-1}_{rr} \xrightarrow{p} \{E(X_r)\}^2 (\mathcal{X}^{-1})_{rr}.$$

Since $\hat{\gamma}_{n0} \xrightarrow{p} \gamma_0$, $\hat{\gamma}_{nr} \xrightarrow{p} \gamma_r$, $s_{\tilde{X}_r}^2 \xrightarrow{p} \text{var}(\tilde{X}_r)$ and $\bar{\tilde{X}}_{nr} \xrightarrow{p} E(X_r)$, the result follows. □

## APPENDIX: AUXILIARY RESULTS AND PROOFS

**A.1. Additional technical conditions.** We introduce some further technical conditions:

(C6) The functions $h_1(u) = \int x g_1(x, u) \, dx$ and $h_2(u) = \int x g_2(x, u) \, dx$ are uniformly Lipschitz, where $g_1(\cdot, \cdot)$ and $g_2(\cdot, \cdot)$ are the joint density functions of $(X, U)$ and $(Xe, U)$, respectively.
(C7) The error term satisfies $E|e^\lambda| < \infty$ for $\lambda > 4$.

Conditions (C1), (C6) and (C7) are needed for the proof of Lemma A.4 given in the next section.

### A.2. Auxiliary results on martingale differences.

LEMMA A.1. *Under the technical conditions* (C1)–(C6), *on event $A_n$* (11) *the martingale differences $Z_{n0t}$ satisfy the conditions*

(a) $$\sum_{t=1}^{n} E\{Z_{n0t}^2 I(|Z_{n0t}| > \varepsilon)\} \to 0 \quad \text{for all } \varepsilon > 0,$$

(b) $$\Delta_{n0}^2 = \sum_{t=1}^{n} Z_{n0t}^2 \xrightarrow{p} \sigma_0^2 \quad \text{for } \sigma_0^2 > 0.$$

PROOF. Let $Z_{n0t} = w_{n0t} v_{n0t}$, where $w_{n0t} = 1/\sqrt{n}$, and

$$v_{n0t} = \gamma_0 \psi(U_{nt}) + \psi(U_{nt}) e_{nt} \{(L_{nj(t)} X'^T_{nj(t)} X'_{nj(t)})^{-1} X'^T_{nj(t)}\}_{1k(t)} - \gamma_0$$

$$= \alpha_{1nt} + \alpha_{2nt} e_{nt},$$

where $\alpha_{1nt} = \gamma_0 \psi(U_{nt}) - \gamma_0$, $\alpha_{2nt} = \psi(U_{nt})\{(L_{nj(t)}^{-1} X'^T_{nj(t)} X'_{nj(t)})^{-1} X'^T_{nj(t)}\}_{1k(t)}$ and $E(v_{n0t}) = 0$. Using (C1), (C3) and (C4), it holds on event $A_n$ that $\sup_{1 \le t \le n} |\alpha_{1nt}| < c_1$ and $\sup_{1 \le t \le n} |\alpha_{2nt}| < c_2$ for some $c_1, c_2 > 0$. Thus, it holds for $\varepsilon > 0$ that

$$\sum_{t=1}^{n} E\{Z_{n0t}^2 I(|Z_{n0t}| > \varepsilon)\} = \sum_{t=1}^{n} \int x^2 I(|x| > \varepsilon) \, dF_{w_{n0t} v_{n0t}}(x)$$

$$= \sum_{t=1}^{n} \int x^2 I(|x| > \varepsilon/|w_{n0t}|) w_{n0t}^2 \, dF_{v_{n0t}}(x)$$

$$= n^{-1} \sum_{t=1}^{n} \int x^2 I(|x| > \sqrt{n}\varepsilon) \, dF_{v_{n0t}}(x)$$

$$\le n^{-1} \sum_{t=1}^{n} \{E(v_{n0t}^4)\}^{1/2} \{P(v_{n0t}^2 > n\varepsilon^2)\}^{1/2}.$$

Now, $E(v_{n0t}^4)$ is bounded uniformly in $n$ and $t$, since $e_{nt}$ has finite fourth moment by (C7), and $P(v_{n0t}^2 > n\varepsilon^2) = P((\alpha_{1nt} + \alpha_{2nt} e_{nt})^2 > n\varepsilon^2) \le P(\alpha_{1nt}^2 + \alpha_{2nt}^2 e_{nt}^2 + 2|\alpha_{1nt} \alpha_{2nt} e_{nt}| > n\varepsilon^2) \le P(c_1^2 + c_2^2 e_{nt}^2 + 2c_1 c_2 |e_{nt}| > n\varepsilon^2)$. Lemma A.1(a) follows, since $P(c_1^2 + c_2^2 e_{nt}^2 + 2c_1 c_2 |e_{nt}| > n\varepsilon^2) \to 0$ uniformly in $n$ and $t$, $e_{nt}^2$ and $|e_{nt}|$ being i.i.d. with finite fourth moments.

The term $\Delta_{n0}^2$ given in Lemma A.1(b) is equal to

$$\Delta_{n0}^2 = \gamma_0^2 \left\{ n^{-1} \sum_t \psi^2(U_{nt}) \right\} + \gamma_0^2 - 2\gamma_0^2 \left\{ n^{-1} \sum_t \psi(U_{nt}) \right\}$$

$$+ 2\gamma_0 n^{-1} \sum_t \psi^2(U_{nt}) e_{nt} \{(L_{nj(t)}^{-1} X'^T_{nj(t)} X'_{nj(t)})^{-1} X'^T_{nj(t)}\}_{1k(t)}$$

$$- 2\gamma_0 n^{-1} \sum_t \psi(U_{nt}) e_{nt} \{(L_{nj(t)}^{-1} X'^T_{nj(t)} X'_{nj(t)})^{-1} X'^T_{nj(t)}\}_{1k(t)}$$

$$+ n^{-1} \sum_t \psi^2(U_{nt}) e_{nt}^2 \{(L_{nj(t)}^{-1} X'^T_{nj(t)} X'_{nj(t)})^{-1} X'^T_{nj(t)}\}_{1k(t)}^2$$

$$= T_1 + \cdots + T_6.$$

It follows from the law of large numbers that

$$T_1 + T_2 + T_3 \xrightarrow{p} \gamma_0^2 E\{\psi^2(U)\} + \gamma_0^2 - 2\gamma_0^2 E\{\psi(U)\} = \gamma_0^2 \operatorname{var}\{\psi(U)\}.$$

On event $A_n$, $E(T_4|U,X,L_{nj}) = 0$ and

$$\operatorname{var}(T_4|U,X,L_{nj}) = \frac{4\sigma^2 \gamma_0^2}{n^2} \sum_t \psi^4(U_{nt}) \{(L_{nj(t)}^{-1} X'^T_{nj(t)} X'_{nj(t)})^{-1} X'^T_{nj(t)}\}_{1k(t)}^2$$

$$= O(n^{-1}).$$

Thus, $E(T_4) = 0$ and $\operatorname{var}(T_4) = O(n^{-1})$, implying that $T_4 = O_p(n^{-1/2})$ on $A_n$. Similarly, it can be shown that $T_5 = O_p(n^{-1/2})$ on $A_n$.

Next consider the last term $T_6$, which can also be written as

$$T_6 = n^{-1} \sum_{j=1}^m \sum_{k=1}^{L_{nj}} \{(L_{nj}^{-1} X'^T_{nj} X'_{nj})^{-1} X'^T_{nj}\}_{1k}^2 \psi^2(U'_{njk}) e'^2_{njk}.$$

Expanding $\{(L_{nj}^{-1} X'^T_{nj} X'_{nj})^{-1} X'^T_{nj}\}_{1k}^2 \psi^2(U'_{njk}) e'^2_{njk}$ for each $k$, we get

$$T_6 = n^{-1} \sum_{j=1}^m \sum_{k=1}^{L_{nj}} \{(L_{nj}^{-1} X'^T_{nj} X'_{nj})^{-1}_{11} e'_{njk} \psi(U'_{njk})$$

$$+ (L_{nj}^{-1} X'^T_{nj} X'_{nj})^{-1}_{12} e'_{njk} \psi(U'_{njk}) X'_{n1jk} + \cdots$$

$$+ (L_{nj}^{-1} X'^T_{nj} X'_{nj})^{-1}_{1,p+1} e'_{njk} \psi(U'_{njk}) X'_{npjk}\}^2,$$

which by Lemma A.4(a) and the law of large numbers is equal to

$$\sigma^2 E\{\psi^2(U)\}[(\mathcal{X}^{-1})^2_{11} + (\mathcal{X}^{-1})^2_{12} E(X_1^2) + \cdots$$

$$+ (\mathcal{X}^{-1})^2_{1,p+1} E(X_p^2)$$

$$+ \{2(\mathcal{X}^{-1})_{11}(\mathcal{X}^{-1})_{12} E(X_1) + \cdots$$

$$+ 2(\mathcal{X}^{-1})_{11}(\mathcal{X}^{-1})_{1,p+1} E(X_p)\}$$

$$+ \{2(\mathcal{X}^{-1})_{12}(\mathcal{X}^{-1})_{13} E(X_1 X_2) + \cdots$$

$$+ 2(\mathcal{X}^{-1})_{12}(\mathcal{X}^{-1})_{1,p+1} E(X_1 X_p)\} + \cdots$$

$$+ \{2(\mathcal{X}^{-1})_{1p}(\mathcal{X}^{-1})_{1,p+1} E(X_{p-1} X_p)\}] + o_p(1)$$

$$= \sigma^2 E\{\psi^2(U)\} (\mathcal{X}^{-1} \mathcal{X}^T \mathcal{X}^{-1^T})_{11} + o_p(1)$$

$$= \sigma^2 E\{\psi^2(U)\} (\mathcal{X}^{-1})_{11} + o_p(1),$$

where $\mathcal{X}$ is as defined in (C5) and given explicitly in (8). Thus

$$\Delta_{n0}^2 \xrightarrow{p} \gamma_0^2 \operatorname{var}\{\psi(U)\} + \sigma^2 (\mathcal{X}^{-1})_{11} E\{\psi^2(U)\} \equiv \sigma_0^2,$$

and Lemma A.1(b) follows. □

LEMMA A.2. *Under the technical conditions* (C1)–(C6), *on event* $A_n$ (11) *the martingale differences* $Z_{nrt}$ *satisfy the conditions*

(a) $\quad \sum_{t=1}^{n} E\{Z_{nrt}^2 I(|Z_{nrt}| > \varepsilon)\} \to 0 \qquad$ *for all* $\varepsilon > 0$,

(b) $\quad \Delta_{nr}^2 = \sum_{t=1}^{n} Z_{nrt}^2 \xrightarrow{p} (a,b)\Sigma_r(a,b)^T \qquad$ *for* $(a,b)\Sigma_r(a,b)^T > 0$.

PROOF. Let $Z_{nrt} = w_{nrt} v_{nrt}$, where $w_{nrt} = 1/\sqrt{n}$, $\alpha_{3nt} = a\gamma_r \psi(U_{nt}) X_{nrt} - a\gamma_r E(X_r) + b\phi_r(U_{nt}) X_{nrt} - bE(X_r)$, $\alpha_{4nt} = a\bar{X}'_{nrj(t)} \psi(U_{nt}) \{(L_{nj(t)}^{-1} X'^T_{nj(t)} \times X'_{nj(t)})^{-1} X'^T_{nj(t)}\}_{rk(t)}$, $v_{nrt} = \alpha_{3nt} + \alpha_{4nt} e_{nt}$ and $E(v_{nrt}) = 0$. On event $A_n$, $\sup_{1 \le t \le n} |\alpha_{3nt}| < c_3$ and $\sup_{1 \le t \le n} |\alpha_{4nt}| < c_4$ for some $c_3, c_4 > 0$, and thus Lemma A.2(a) follows in a fashion similar to Lemma A.1(a).

The term $\Delta_{nr}^2$ in Lemma A.2(b) is equal to

$$\Delta_{nr}^2 = a^2 \gamma_r^2 \left\{ n^{-1} \sum_t \psi^2(U_{nt}) X_{nrt}^2 \right\} + a^2 \gamma_r^2 \{E(X_r)\}^2 + b^2 \left\{ n^{-1} \sum_t \phi_r^2(U_{nt}) X_{nrt}^2 \right\}$$

$$+ b^2 \{E(X_r)\}^2 - 2a^2 \gamma_r^2 E(X_r) \left\{ n^{-1} \sum_t \psi(U_{nt}) X_{nrt} \right\} + 2ab\gamma_r \{E(X_r)\}^2$$

$$+ 2ab\gamma_r \left\{ n^{-1} \sum_t \psi(U_{nt}) \phi_r(U_{nt}) X_{nrt}^2 \right\} - 2b^2 E(X_r) \left\{ n^{-1} \sum_t \phi_r(U_{nt}) X_{nrt} \right\}$$

$$- 2ab\gamma_r E(X_r) \left\{ n^{-1} \sum_t \psi(U_{nt}) X_{nrt} \right\}$$

$$- 2ab\gamma_r E(X_r) \left\{ n^{-1} \sum_t \phi_r(U_{nt}) X_{nrt} \right\}$$

$$+ 2a^2 \gamma_r n^{-1} \sum_t \psi^2(U_{nt}) e_{nt} \bar{X}'_{nrj(t)} X_{nrt} \{(L_{nj(t)}^{-1} X'^T_{nj(t)} X'_{nj(t)})^{-1} X'^T_{nj(t)}\}_{rk(t)}$$

$$- 2a^2 \gamma_r E(X_r) n^{-1} \sum_t \psi(U_{nt}) e_{nt} \bar{X}'_{nrj(t)} \{(L_{nj(t)}^{-1} X'^T_{nj(t)} X'_{nj(t)})^{-1} X'^T_{nj(t)}\}_{rk(t)}$$

$$+ 2abn^{-1} \sum_t \psi(U_{nt}) \phi_r(U_{nt}) e_{nt} \bar{X}'_{nrj(t)} X_{nrt} \{(L_{nj(t)}^{-1} X'^T_{nj(t)} X'_{nj(t)})^{-1} X'^T_{nj(t)}\}_{rk(t)}$$

$$- 2abE(X_r)n^{-1}\sum_t \psi(U_{nt})e_{nt}\bar{X}'_{nrj(t)}\{(L^{-1}_{nj(t)}X'^T_{nj(t)}X'_{nj(t)})^{-1}X'^T_{nj(t)}\}_{rk(t)}$$

$$+ a^2n^{-1}\sum_t \psi^2(U_{nt})e^2_{nt}\bar{X}'^2_{nrj(t)}\{(L^{-1}_{nj(t)}X'^T_{nj(t)}X'_{nj(t)})^{-1}X'^T_{nj(t)}\}^2_{rk(t)}$$

$$= T_1 + \cdots + T_{15},$$

and by the law of large numbers

$$T_1 + \cdots + T_{10} \xrightarrow{p} a^2\gamma^2_r[\{E(X_r)\}^2 \mathrm{var}\{\psi(U)\} + \mathrm{var}(X_r)E\{\psi^2(U)\}]$$
$$+ 2ab\gamma_r[E\{\phi_r(U)\psi(U)\}E(X^2_r) - \{E(X_r)\}^2] + b^2\mathrm{var}(\tilde{X}_r).$$

On event $A_n$, $E(T_{11}|U,X,L_{nj}) = 0$ and

$$\mathrm{var}(T_{11}|U,X,L_{nj})$$
$$= \frac{4a^4\sigma^2\gamma^2_r}{n^2}\sum_t \psi^4(U_{nt})\bar{X}'^2_{nrj}X'^2_{nrt}\{(L^{-1}_{nj(t)}X'^T_{nj(t)}X'_{nj(t)})^{-1}X'^T_{nj(t)}\}^2_{rk(t)},$$

which is $O(n^{-1})$. Thus, $E(T_{11}) = 0$ and $\mathrm{var}(T_{11}) = O(n^{-1})$, implying that $T_{11} = O_p(n^{-1/2})$ on $A_n$. Similarly, it can be shown that $T_{12} = T_{13} = T_{14} = O_p(n^{-1/2})$ on $A_n$.

Next consider the last term $T_{15}$, which can also be expressed as

$$T_{15} = a^2n^{-1}\sum_{j=1}^m\sum_{k=1}^{L_{nj}}\{(L^{-1}_{nj}X'^T_{nj}X'_{nj})^{-1}X'^T_{nj}\}^2_{rk}\psi^2(U'_{njk})e'^2_{njk}\bar{X}'^2_{nrj}.$$

Again expanding $\{(L^{-1}_{nj}X'^T_{nj}X'_{nj})^{-1}X'^T_{nj}\}^2_{1k}\psi^2(U'_{njk})e'^2_{njk}\bar{X}'^2_{nrj}$ for each $k$, we get

$$T_{15} = a^2n^{-1}\sum_{j=1}^m\sum_{k=1}^{L_{nj}}\{(L^{-1}_{nj}X'^T_{nj}X'_{nj})^{-1}_{r1}\bar{X}'_{nrj}e'_{njk}\psi(U'_{njk})$$
$$+ (L^{-1}_{nj}X'^T_{nj}X'_{nj})^{-1}_{r2}\bar{X}'_{nrj}e'_{njk}\psi(U'_{njk})X'_{n1jk} + \cdots$$
$$+ (L^{-1}_{nj}X'^T_{nj}X'_{nj})^{-1}_{r,p+1}\bar{X}'_{nrj}e'_{njk}\psi(U'_{njk})X'_{npjk}\}^2,$$

which by Lemma A.4(a) and the law of large numbers is equal to

$$a^2\sigma^2\{E(X_r)\}^2E\{\psi^2(U)\}$$
$$\times [(\mathcal{X}^{-1})^2_{r1} + (\mathcal{X}^{-1})^2_{r2}E(X^2_1) + \cdots + (\mathcal{X}^{-1})^2_{r,p+1}E(X^2_p)$$
$$+ \{2(\mathcal{X}^{-1})_{r1}(\mathcal{X}^{-1})_{r2}E(X_1) + \cdots + 2(\mathcal{X}^{-1})_{r1}(\mathcal{X}^{-1})_{r,p+1}E(X_p)\}$$
$$+ \{2(\mathcal{X}^{-1})_{r2}(\mathcal{X}^{-1})_{r3}E(X_1X_2) + \cdots$$
$$+ 2(\mathcal{X}^{-1})_{r2}(\mathcal{X}^{-1})_{r,p+1}E(X_1X_p)\} + \cdots$$

$$+ \{2(\mathcal{X}^{-1})_{rp}(\mathcal{X}^{-1})_{r,p+1}E(X_{p-1}X_p)\}]$$
$$+ o_p(1)$$
$$= a^2\sigma^2\{E(X_r)\}^2 E\{\psi^2(U)\}(\mathcal{X}^{-1}\mathcal{X}^T\mathcal{X}^{-1^T})_{rr} + o_p(1)$$
$$= a^2\sigma^2\{E(X_r)\}^2 E\{\psi^2(U)\}(\mathcal{X}^{-1})_{rr} + o_p(1).$$

Thus
$$\Delta^2_{nr} \xrightarrow{p} (a,b)\Sigma_r(a,b)^T = (a,b)\begin{bmatrix} \Sigma_{r11} & \Sigma_{r12} \\ \Sigma_{r12} & \Sigma_{r22} \end{bmatrix}(a,b)^T,$$

where $\Sigma_{r11} = \gamma_r^2[\{E(X_r)\}^2 \operatorname{var}\{\psi(U)\} + \operatorname{var}(X_r)E\{\psi^2(U)\}] + \sigma^2\{E(X_r)\}^2 \times E\{\psi^2(U)\}(\mathcal{X}^{-1})_{rr}$, $\Sigma_{r12} = \gamma_r[E\{\phi_r(U)\psi(U)\}E(X_r^2) - \{E(X_r)\}^2]$ and $\Sigma_{r22} = \operatorname{var}(\tilde{X}_r)$. Hence Lemma A.2(b) follows. □

**A.3. Auxiliary results on approximations of inverses.** Defining $\delta_{n0jk} = \psi(U'_{njk}) - \psi(U'^*_{nj})$ and $\delta_{nrjk} = \phi_r(U'_{njk}) - \phi_r(U'^*_{nj})$ for $1 \le k \le L_j$ and $1 \le r \le p$, where $U'^*_{nj} = L_{nj}^{-1}\sum_{k=1}^{L_{nj}} U'_{njk}$ is the average of the $U$'s in $B_{nj}$, we obtain the following results, by Taylor expansions and boundedness considerations: for $1 \le t, s \le p$, $0 \le r, r' \le p$ and $1 \le \ell \le 2$, (a) $\sup_{k,j}|U'_{njk} - U'^*_{nj}| \le (b-a)/m$; (b) $\sup_{k,j}|\delta_{nrjk}| = O(m^{-1})$; (c) $\sup_j|L_{nj}^{-1}\sum_k \delta_{nrjk} X'^\ell_{ntjk}| = O(m^{-1})$; (d) $\sup_j|L_{nj}^{-1}\sum_k \delta^2_{nrjk} X'^\ell_{ntjk}| = O(m^{-2})$; (e) $\sup_j|L_{nj}^{-1}\sum_k \delta_{nrjk} X'_{ntjk} X'_{nsjk}| = O(m^{-1})$; (f) $\sup_j|L_{nj}^{-1}\sum_k \delta_{nrjk} \times \delta_{nr'jk} X'_{ntjk} X'_{nsjk}| = O(m^{-2})$.

LEMMA A.3. *Under the technical conditions* (C1)–(C6), *it holds on event* $E_n$ (12) *that*
$$\sup_j |(L_{nj}^{-1}\tilde{X}'^T_{nj}\tilde{X}'_{nj})^{-1} - (\Phi_{nj} \boxdot \Xi_{nj})| = O(m^{-1})\mathbf{1}_{(p+1)\times(p+1)},$$

*where*

$$(23) \quad \Phi_{nj} = \begin{bmatrix} 1 & 1/\phi_1(U'^*_{nj}) & \cdots & 1/\phi_p(U'^*_{nj}) \\ 1/\phi_1(U'^*_{nj}) & 1/\phi_1^2(U'^*_{nj}) & \cdots & 1/(\phi_p(U'^*_{nj})\phi_1(U'^*_{nj})) \\ \vdots & & \ddots & \\ 1/\phi_p(U'^*_{nj}) & 1/(\phi_p(U'^*_{nj})\phi_1(U'^*_{nj})) & \cdots & 1/\phi_p^2(U'^*_{nj}) \end{bmatrix},$$

$\Xi_{nj} = (L_{nj}^{-1} X'^T_{nj} X'_{nj})^{-1}$ *and* $\mathbf{1}_{(p+1)\times(p+1)}$ *denotes the* $(p+1)\times(p+1)$ *matrix of 1's.*

PROOF. The proof is by induction on $p$. Define

$$(24) \quad \tilde{X}'^{(\ell)}_{nrj} = \frac{1}{L_{nj}}\sum_{k=1}^{L_{nj}} \tilde{X}'^\ell_{nrjk}, \qquad (\tilde{X}'_{nrj}\tilde{X}'_{nsj})^{(\ell)} = \frac{1}{L_{nj}}\sum_{k=1}^{L_{nj}} (\tilde{X}'_{nrjk}\tilde{X}'_{nsjk})^\ell,$$

and analogously for $X'^{(\ell)}_{nrj}$ and $(X'_{nrj}X'_{nsj})^{(\ell)}$ where $1 \le r, s \le p$. First consider the claim for $p = 1$ on $E_n$,

$$(L_{nj}^{-1}\tilde{X}'^T_{nj}\tilde{X}'_{nj})^{-1} = \frac{1}{\tilde{X}'^{(2)}_{n1j} - (\tilde{X}'^{(1)}_{n1j})^2} \begin{bmatrix} \tilde{X}'^{(2)}_{n1j} & -\tilde{X}'^{(1)}_{n1j} \\ -\tilde{X}'^{(1)}_{n1j} & 1 \end{bmatrix}.$$

By boundedness considerations and properties (c) and (d), it holds that $\sup_j |\tilde{X}'^{(2)}_{n1j} - \phi_1^2(U'^*_{nj})X'^{(2)}_{n1j}| = O(m^{-1})$, $\sup_j |\tilde{X}'^{(1)}_{n1j} - \phi_1(U'^*_{nj})X'^{(1)}_{n1j}| = O(m^{-1})$, and therefore

$$\sup_j |\{\tilde{X}'^{(2)}_{n1j} - (\tilde{X}'^{(1)}_{n1j})^2\} - \phi_1^2(U'^*_{nj})\{X'^{(2)}_{n1j} - (X'^{(1)}_{n1j})^2\}| = \sup_j |\tilde{d}_{nj} - \phi_1^2(U'^*_{nj})d_{nj}|$$

$$= O(m^{-1}),$$

where $\tilde{d}_{nj} = \det(L_{nj}^{-1}\tilde{X}'^T_{nj}\tilde{X}'_{nj})$ and $d_{nj} = \det(L_{nj}^{-1}X'^T_{nj}X'_{nj})$. Thus,

$$\sup_j |(L_{nj}^{-1}\tilde{X}'^T_{nj}\tilde{X}'_{nj})^{-1} - (\Phi_{nj} \boxdot \Xi_{nj})| = O(m^{-1})\mathbf{1}_{2\times 2},$$

where $(\Phi_{nj})_{2\times 2}$ is as given in (23) and $(\Xi_{nj})_{2\times 2} = (L_{nj}^{-1}X'^T_{nj}X_{nj})^{-1}_{2\times 2}$.

Next, we show that Lemma A.3 holds for $p+1$, assuming it holds for $p$. Let

$$(L_{nj}^{-1}\tilde{X}'^T_{nj}\tilde{X}'_{nj})_{(p+2)\times(p+2)} = B_{nj} = \begin{bmatrix} B_{nj11} & B_{nj12} \\ B^T_{nj12} & B_{nj22} \end{bmatrix},$$

$$(L_{nj}^{-1}\tilde{X}'^T_{nj}\tilde{X}'_{nj})^{-1}_{(p+2)\times(p+2)} = B_{nj}^{-1} = \begin{bmatrix} B^{11}_{nj} & B^{12}_{nj} \\ B^{12^T}_{nj} & B^{22}_{nj} \end{bmatrix},$$

and similarly let

$$(L_{nj}^{-1}X'^T_{nj}X'_{nj})_{(p+2)\times(p+2)} = D_{nj} = \begin{bmatrix} D_{nj11} & D_{nj12} \\ D^T_{nj12} & D_{nj22} \end{bmatrix},$$

$$(L_{nj}^{-1}X'^T_{nj}X'_{nj})^{-1}_{(p+2)\times(p+2)} = D_{nj}^{-1} = \begin{bmatrix} D^{11}_{nj} & D^{12}_{nj} \\ D^{12^T}_{nj} & D^{22}_{nj} \end{bmatrix},$$

where $B_{nj11} = (L_{nj}^{-1}\tilde{X}'^T_{nj}\tilde{X}'_{nj})_{(p+1)\times(p+1)}$ and $D_{nj11} = (L_{nj}^{-1}X'^T_{nj}X'_{nj})_{(p+1)\times(p+1)}$. By the assumption,

(25) $\quad \sup_j |B^{-1}_{nj11} - (\Phi_{nj} \boxdot \Xi_{nj})_{(p+1)\times(p+1)}| = O(m^{-1})\mathbf{1}_{(p+1)\times(p+1)}.$

By properties (c), (d), (e), (f) and boundedness considerations, it holds that

(26) $\quad \sup_j |B_{nj12} - (V_{nj} \boxdot D_{nj12})| = O(m^{-1})\mathbf{1}_{(p+1)\times 1},$

(27) $\quad \sup_j |B_{nj22} - \phi^2_{n(p+1)}(U'^*_{nj})D_{nj22}| = O(m^{-1}),$

where $B_{nj12}^T = (\tilde{X}'^{(1)}_{n(p+1)j}, (\tilde{X}'_{n(p+1)j}\tilde{X}'_{n1j})^{(1)}, \ldots, (\tilde{X}'_{n(p+1)j}\tilde{X}'_{npj})^{(1)})$, $D_{nj12}^T = (X'^{(1)}_{n(p+1)j}, (X'_{n(p+1)j}X'_{n1j})^{(1)}, \ldots, (X'_{n(p+1)j}X'_{npj})^{(1)})$, $B_{nj22} = \tilde{X}'^{(2)}_{n(p+1)j}$, $D_{nj22} = X'^{(2)}_{n(p+1)j}$ and $V_{nj}^T = (\phi_{p+1}(U'^*_{nj}), \phi_{p+1}(U'^*_{nj})\phi_1(U'^*_{nj}), \ldots, \phi_{p+1}(U'^*_{nj})\phi_p(U'^*_{nj}))$. Since $B_{nj}^{22} = (B_{nj22} - B_{nj12}^T B_{nj11}^{-1} B_{nj12})^{-1}$, using (25), (26), (27) and the uniform boundedness of $D_{nj12}$, $D_{nj11}^{-1}$, $D_{nj22}$ on $A_n$,

$$\sup_j |B_{nj}^{22} - \{\phi_{p+1}^2(U'^*_{nj})D_{nj}^{22}\}^{-1}| = O(m^{-1}),$$

where $\inf_j |\phi_{p+1}^2(U'^*_{nj})D_{nj}^{22}| = \inf_j |\phi_{p+1}^2(U'^*_{nj})(D_{nj22} - D_{nj12}^T D_{nj11}^{-1} D_{nj12})| > 0$, since $\phi_{p+1}(\cdot)$ is assumed to be strictly positive, and since $\sup_j |D_{nj22} - D_{nj12}^T \times D_{nj11}^{-1} D_{nj12}| > 0$. The latter holds on $A_n$, since then $\sup_j |d_j| = \sup_j |\det(D_{nj11}) \times (D_{nj22} - D_{nj12}^T D_{nj11}^{-1} D_{nj12})| > 0$.

Now $B_{nj}^{11} = B_{nj11}^{-1} + B_{nj11}^{-1} B_{nj12} B_{nj}^{22} B_{nj12}^T B_{nj11}^{-1}$. Since $D_{nj12}$, $D_{nj11}^{-1}$ are uniformly bounded on $A_n$,

(28) $$\sup_j |B_{nj}^{11} - (\Phi_{nj} \boxdot \Gamma_{nj})| = O(m^{-1}) \mathbf{1}_{(p+1) \times (p+1)},$$

where $\Phi_{nj}$ is as defined in (23), and $\Gamma_{nj} = D_{nj11}^{-1} + D_{nj11}^{-1} D_{nj12} D_{nj}^{22} D_{nj12}^T D_{nj11}^{-1} = D_{nj}^{11}$.

Since $B_{nj}^{12} = -B_{nj}^{11} B_{nj12} B_{nj22}^{-1}$, using (26), (27), (28) and boundedness considerations,

$$\sup_j |B_{nj}^{12} - (\Omega_{nj} \boxdot \Lambda_{nj})| = O(m^{-1}) \mathbf{1}_{(p+1) \times 1},$$

where $\Omega_{nj}^T = (1/\phi_{p+1}(U'^*_{nj}), 1/\{\phi_{p+1}(U'^*_{nj})\phi_1(U'^*_{nj})\}, \ldots, 1/\{\phi_{p+1}(U'^*_{nj})\phi_p(U'^*_{nj})\})$ and $\Lambda_{nj} = -D_{nj}^{11} D_{nj12} D_{nj22}^{-1} = D_{nj}^{12}$. Thus, reassembling the partitioned matrix $B_{nj}^{-1}$, Lemma A.3 follows. □

LEMMA A.4. *Under the technical conditions* (C1)–(C7), *for a sequence $r_n$ such that $r_n = O_p\{\sqrt{(m\log n)/n}\}$, on event $A_n$ (11)*

(a) $$\sup_j |(L_{nj}^{-1} X'^T_{nj} X'_{nj})^{-1} - \mathcal{X}^{-1}| = O_p(r_n) \mathbf{1}_{(p+1) \times (p+1)},$$

(b) $$\sup_j |L_{nj}^{-1} X'^T_{nj} e'_{nj}| = O_p(r_n) \mathbf{1}_{(p+1) \times 1},$$

where $\mathcal{X}$ as defined in (8) is assumed to be nonsingular by (C5), and $e'_{nj} = (e'_{nj1}, \ldots, e'_{njL_{nj}})^T$.

PROOF. Using the sample moment notation in (24),

$$\frac{1}{L_{nj}}X_{nj}'^T X_{nj}' = \begin{bmatrix} 1 & X_{n1j}'^{(1)} & \cdots & X_{npj}'^{(1)} \\ X_{n1j}'^{(1)} & X_{n1j}'^{(2)} & \cdots & (X_{n1j}' X_{npj}')^{(1)} \\ \vdots & & \ddots & \\ X_{npj}'^{(1)} & (X_{npj}' X_{n1j}')^{(1)} & \cdots & X_{npj}'^{(2)} \end{bmatrix}_{(p+1)\times(p+1)}$$

leads to

$$d_{nj} = \sum (-1)^{\text{sign}(\tau)} (L_{nj}^{-1} X_{nj}'^T X_{nj}')_{1\tau(1)} \cdots (L_{nj}^{-1} X_{nj}'^T X_{nj}')_{(p+1),\tau(p+1)},$$

where the sum is taken over all permutations $\tau$ of $(1,\ldots,p+1)$, and $\text{sign}(\tau)$ equals $+1$ or $-1$, depending on whether $\tau$ can be written as the product of an even or odd number of transpositions. The terms in the above sum have the general form

$$(29) \qquad X_{nr_1j}'^{(1)}(X_{n1j}' X_{nr_2j}')^{(1)} \cdots (X_{npj}' X_{nr_{p+1}j}')^{(1)},$$

where $X_0' = 1$ and $(r_1,\ldots,r_{p+1})$ is a permutation of $(0,\ldots,p)$. Considering the definition of the Nadaraya–Watson kernel estimator [10, 12], we note that an arbitrary term in (29) has the form $(X_{nsj}' X_{nr_{s+1}j}')^{(1)} = \hat{m}_{nsr_{s+1}}(U_{nj}^M)$ for $0 \le s \le p+1$, $K(\cdot) = (1/2)\mathbf{1}_{[-1,1]}$, $h = (b-a)/m$, and $U_{nj}^M = a + (2j-1)\{(b-a)/(2m)\}$ are the midpoints of the bins $B_{nj}$. Uniform consistency of Nadaraya–Watson estimators with kernels of compact support has been shown in [4], where

$$(30) \qquad \sup_{a \le u \le b} |\hat{m}_{nsr_{s+1}}(u) - m_{sr_{s+1}}(u)| = O_p(r_n),$$

$m_{sr_{ns+1}}(u) = E(X_s X_{r_{s+1}} | U = u) = E(X_s X_{r_{s+1}})$, and $r_n$ is as defined in Lemma A.4. Then (30) implies

$$(31) \qquad \begin{aligned} \sup_j |\hat{m}_{nsr_{s+1}}(U_{nj}^M) - m_{sr_{s+1}}(U_{nj}^M)| &= O_p(r_n), \\ \sup_j |(X_{nsj}' X_{nr_{s+1}j}')^{(1)} - E(X_s X_{r_{s+1}})| &= O_p(r_n). \end{aligned}$$

Hence the uniform consistency of (29) follows, where the limit of (29) is $E(X_{r_1})E(X_1 X_{r_2})\cdots E(X_p X_{r_{p+1}})$, and

$$(32) \qquad \sup_j |d_{nj} - \det(\mathcal{X})| = O_p(r_n)$$

follows.

The cofactor of $(L_{nj}^{-1} X_{nj}'^T X_{nj})_{r\ell}$ is defined by $(-1)^{r+\ell}$ times the minor of $(L_{nj}^{-1} X_{nj}'^T X_{nj})_{r\ell}$, where the minor is the determinant after deleting the $r$th row and the $\ell$th column of $(L_{nj}^{-1} X_{nj}'^T X_{nj})$. With a similar argument as in

INFERENCE FOR ADJUSTED REGRESSION 25the case of $d_{nj}$, it can be shown that the minor of $(L_{nj}^{-1} X_{nj}'^T X_{nj})_{r\ell}$ converges uniformly over $j$ to the minor of $(\mathcal{X})_{r\ell}$ with rate $r_n$. Thus part (a) of the lemma follows. For part (b) of the lemma, consider

$$L_{nj}^{-1} X_{nj}'^T e_{nj}' = \left( L_{nj}^{-1} \sum_{k=1}^{L_{nj}} e_{njk}', L_{nj}^{-1} \sum_{k=1}^{L_{nj}} X_{n1jk}' e_{njk}', \ldots, L_{nj}^{-1} \sum_{k=1}^{L_{nj}} X_{npjk}' e_{njk}' \right)^T.$$

Each term in the above sum is equal to $\hat{m}(U_{nj}^M)$, where $m(U_{nj}^M) = E(e|U) = 0$ or $m(U_{nj}^M) = E(X_r e|U) = 0$, for $r = 1, \ldots, p$. Thus by the uniform consistency of $\hat{m}(U_{nj}^M)$ for $m(U_{nj}^M)$, part (b) of the lemma follows.

On event $A_n$, (32) implies that $P(\inf_j d_{nj} > \zeta) \to 1$ as $n \to \infty$, where $\zeta = \min\{\rho/2, [\inf_j(\phi_1^2(U_{nj}'^*), \ldots, \phi_p^2(U_{nj}'^*))]^p \rho/2\}$ and $\rho$ is as defined in (C5). We also need to show $P(\min_j L_{nj} \leq p) \to 0$ as $n \to \infty$ in order to show that $P(A) \to 1$ as $n \to \infty$. Since $P(\min_j L_{nj} > p) = 1 - P(0 \leq L_{nj} \leq p$ for all $j = 1, \ldots, m) \geq 1 - \sum_{j=1}^{m} P(0 \leq L_{nj} \leq p) \geq 1 - m \sup_j P(0 \leq L_{nj} \leq p)$, it is enough to show $P(0 \leq L_{nj} \leq p) = o(m^{-1})$ uniformly in $j$. Now, $L_{nj} \sim \text{Bin}(n, p_{nj})$, where $c_1(b-1)/m \leq p_{nj} \leq c_2(b-a)/m$ uniformly in $j$, and $c_1, c_2$ are as given in (C1). Therefore, $mP(0 \leq L_{nj} \leq p) = m \sum_{x=0}^{p} p_{nj}^x (1 - p_{nj})^{n-x} n!/(x!(n-x)!) \leq m \sum_{x=0}^{p} n^x \{c_2(b-a)/m\}^x \{1-(c_1(b-a)/m)\}^{n-x} \approx \sum_{x=0}^{p} m(n/m)^x \{e^{-c_1(b-a)}\}^{n/m}$, where "≈" is used to denote asymptotic equivalence. The previously made assumption of $m \log n/n \to 0$ as $n \to 0$ implies $\log m/(n/m) \to 0$ as $n \to 0$. Thus, $\log m + x \log(n/m) - nc_1(b-a)/m \to -\infty$, $m(n/m)^x \{e^{-c_1(b-a)}\}^{n/m} \to 0$ for $x = 0, \ldots, p$ and $mP(0 \leq L_{nj} \leq p) \to 0$ uniformly in $j$ as $n \to \infty$. It follows that $P(A) \to 1$ as $n \to \infty$.

Furthermore, Lemma A.3 implies

$$\sup_j |\tilde{d}_{nj} - \phi_1^2(U_{nj}'^*) \cdots \phi_p^2(U_{nj}'^*) d_{nj}| = O_p(m^{-1}).$$

This shows that $P(\inf_j \tilde{d}_{nj} > \zeta) \to 1$ as $n \to \infty$, which implies $P(\tilde{A}_n) \to 1$ as $n \to \infty$. Thus $P(E_n) \to 1$ as $n \to \infty$. $\square$

Acknowledgments**Acknowledgments.** We wish to thank an Associate Editor and two referees for helpful comments.

DEPARTMENT OF STATISTICS
PENNSYLVANIA STATE UNIVERSITY
411 JOAB L. THOMAS BUILDING
UNIVERSITY PARK, PENNSYLVANIA 16802
USA
E-MAIL: dsenturk@stat.psu.edu

DEPARTMENT OF STATISTICS
UNIVERSITY OF CALIFORNIA, DAVIS
ONE SHIELDS AVENUE
DAVIS, CALIFORNIA 95616
USA
E-MAIL: mueller@wald.ucdavis.edu